# EPISTEMOGRAPHY AND ALGEBRA

Jean-Philippe Drouhard[1]

University of Nice

(IREM de Nice, IUFM Célestin Freinet, UMR P3 "ADEF")

*We propose to address the problem of how to know students' knowledge in an entirely new approach called "epistemography" which is, roughly, an attempt to describe the structure of this knowledge. We claim that what is to be known is made of five tightly interrelated organised systems: the mathematical universe, the system of semio-linguistic representations, the instruments, the rules of the mathematical game, and the identifiers.*

One of the most commonly shared principle of didactics of mathematics is that teaching must ground on students' previous knowledge. Therefore we researchers (and teachers too!) need to know what students know and what they are supposed to know.

But the point is that knowing what students are supposed to know is less easy to do than it appears at a first glance, particularly when they shift from primary studies to secondary studies and when there are frequent curricular changes in the primary studies. In this case, secondary teachers cannot rely on remembering their primary school time; reading curricular documents is not very helpful, neither discussing with primary teachers. The problem is the lack of a common language, or better said, that the common language is not accurate enough. Saying that "students know the sense of operations" or that they are able to solve "simple word problems" is far too fuzzy and superficial.

We propose to address this problem (how to know students' knowledge) in an entirely new approach called "epistemography" which is, roughly, an attempt to describe the structure of this knowledge.

Epistemography is based on an attempt to generalise and conceptualise findings about knowledge we made mainly during previous researches on algebraic thinking. According with many authors we found that semiotic and linguistic knowledge plays a central role in Algebraic Thinking. And we faced the following question: to what extent is this knowledge, *mathematical*? Letters and symbols are not mathematical objects in the same way that numbers or sets or functions are[2]; but on the other hand

1 JPDrouhard@gmail.com

2 More precisely, digits, letters, symbols and expressions made with them form a "language". Languages are they are equally necessary to do mathematics.

Epistemography is a description of the structure of what the subjects have to know in order to actually do mathematics (and not just to pretend to do mathematics!). We chose to call this theory "epistemography" because it is about knowledge ("epistemo-") but, unlike epistemology, not in a historical perspective: rather, epistemography is a kind of geography of knowledge.

We claim that what is to be known is made of five tightly interrelated organised systems: the mathematical universe, the system of semio-linguistic representations, the instruments, the rules of the mathematical game, and the identifiers. We will now present in detail these five knowledge systems. Due to the lack of space this presentation is a quite schematic and abstract one; a much more detailed and discussed presentation of epistemography is to be written.

## THE MATHEMATICAL UNIVERSE

To solve some algebraic problems, you must know that the product of two negative numbers is positive. You can believe that negative numbers are real numbers, or just "imaginary" ones; whatever philosophical option you take, if you want to do mathematics, you need to have some knowledge about something. We call a "mathematical object" this "something", and the Mathematical Universe the system made up of these mathematical objects (e.g. numbers), their relations (e. g. rational numbers are real numbers) and properties (e. g. the product of two negative numbers is positive). Usually, objects of the mathematical universe may be described as individuals (like the number 20) or classes (the even numbers).

## SEMIO-LINGUISTIC REPRESENTATIONS SYSTEM

How to avoid, however, considering as belonging to the mathematical universe, objets or properties whose nature is totally different? We must, actually, distinguish carefully (mathematical) objects (like the number 20) from their (semiolinguistic[3]) representations (like the string of characters "20" made of a "2" and a "0", but also "XX" made of two "X" or ":::::  :::::" made of twenty dots). This distinction –and its consequences– is essential and has been stressed by many authors (Drouhard & Teppo, 2004, Duval, 1995, 2000, 2006, Ernest, 2006, Kirshner, 1989, Radford, 2006, Bagni, 2007 amongst many others). Misunderstanding or neglecting this distinction may lead to quite severe consequences on mathematics learning and teaching studies. Hence our claim is that, besides knowledge about objects of mathematical universe, students must have some (at least practical) knowledge of the very complex and heterogeneous, and often hidden, system of semio-linguistic representations.

mathematically described by the "Language Theory" (a part of Mathematical Logic, shared with computer science).

[3] "semio-" means "related to signs" and "linguistic", "related to language"; see further.

But, how can we decide if a given property is mathematical or semio-linguistic? There is a practical criterion: mathematic properties may be called "representation-free": they remain true whatever representation system is used. For example, the irrationality of $\sqrt{2}$ does not depend on how integers, square roots or fractions are written. Actually the Greeks' notations of the first proof had nothing in common with ours (in particular they did not use any symbolic writing). Semiotic properties, on the contrary, rely on representational conventions. The property that in order to write 1/3 you need an infinite number of (decimal) digits is true – in base ten only; it is false in base three ("0,1": zero unit and one third) or, as in the Babylonian system, in base sixty:𒌋𒌋: two times ten sixtieths.

## Mathematical language

What are the characteristics of the semio-linguistic system? First of all, the "mathematical language" (in a loose sense) is a written one[4]. Mathematical semio-linguistic *units* are written texts. Following and extending Laborde's ideas (1990), written mathematical texts are heterogeneous, made of natural language sentences, symbolic writings, diagrams and tables, graphs and illustrations. Their organisation follows what we call the fruit cake analogy, the natural language being the dough and the symbolic writings, diagrams, graphs and illustrations being the fruit pieces. To describe rigorously such a complex structure is far from easy.

## Linguistic system

Students' ability to understand natural language mathematical texts (the "dough") is linguistic by nature. Mathematical natural language (we call it the "mathematicians jargon") is mostly the natural language itself; but Laborde (1982) showed there are some differences (unusual syntactic constructions like "Let *x* be a number...") between the jargon and the mother-tongue, difficult to interpret by students.

Symbolic writings (like "$b^2 - 4ac > 0$") make up a language, too (Brown & Drouhard, 2004, Drouhard et al, 2006), which is far more complex and different from mother-tongue than it appears at first sight; detailed and accurate descriptions of this language can be found in Kirshner (1987) and Drouhard (1992). Students must learn this language and its syntax[5] – which allows symbolic manipulation (Bell, 1996): the actual mathematic language, ruled by a rigid syntax, permits to perform operations on the symbolic expressions rather than on (mental or graphic) representations.

The present mathematical language is also characterised by a complex but precise

[4] which puts upside down the usual relationship between oral speech and written texts

[5] the syntax is the part of the grammar which deals with the rules that relate one to another the elements of a language. (Syntax says that a parenthesis must be close once opened...

semantics. Semantics (the science of the meaning) is the set of rules and procedures which allows interpreting expressions, in other words which allows relating expressions to mathematical objects.

The most accurate description of this semantics (how symbolic writings refer to mathematical objects and properties) is based on G. Frege's ideas (Drouhard, 1995). G. Frege's key concepts are "denotation" (which can be a numerical value (in the case of "20"), a numerical function (in the case of "$x$+1"), a truth value (in the case of "1 > 20"□) or a boolean function ("$x$+1> 20"□), according to the type of symbolic writing)[6], and "sense□" (the way denotation is given). The linguistic nature of students□□' difficulties with symbolic writings is often underestimated, or confused with conceptual difficulties.

### Semiotic system

Let's give an example of a semiotic problem in algebra. How to represent an infinite series of decimals? Imagine I ask you what the properties of the number 0,666… are. When multiplied by 3 it gives 2? No. Actually I had in mind the number 1999/3000. And yes, I cheated: I broke the representational rule of decimals, which is a semiotic rule (on how to interpret elements like "…") about linguistic objects (the numeric expressions).

There are more than one approach to mathematics semiotics, which were fully presented in the special issue N° 134 (2003) of *Educational Studies in Mathematics*. Duval dedicated his lifelong work to an extensive and coherent theory of semiotics of mathematics education. Three key concepts are the semiotic representation registers, the treatments (within a register) and the conversions (between different registers). Other researchers (see amongst others Otte, 2006) are investigating how to interpret mathematics education using the terms of the founder of semiotics, Charles S. Peirce (1991): the three types of signs –index, icon, symbol– and, maybe more interesting, the three types of inferences –induction, abduction, deduction).

An entire communication paper would not suffice to present even a small part of the outcomes of semiotics for the study of algebraic thinking. Hence we called "semio-linguistic" the mathematics representation system. Therefore students must handle both aspects of this representation system, the linguistic as well as the semiotic one, and the complex interaction between them.

## INSTRUMENTS

Up to now we have seen that to do mathematics, students must not only know objects and how to represent them: now we will see that they need also to know how to use instruments (Rabardel & Vérillon, 1995) to operate on the representations of objects.

[6] The AlNuSet software, developed by Giampaolo Chiappini allows (in a totally original way) a dynamic view of the denotation of algebraic expressions.

However, unlike object/representation opposition, instruments are not characterised by their nature (mathematical objects can also be tools, as noted by Douady, 1986) but instead by their use. Students, then, must learn what these instruments are and how to use them. Given that instruments are only characterised by their use, it is possible to propose a typology, based on their nature: material instruments (like rulers or compasses, see Bagni, 2007), conceptual instruments (mathematical properties, like theorems), semiotic instruments (manipulations on semiotic representations) – this idea appears in L. S. Vygotsky, 1986); eventually one may consider “meta” instruments like strategies and, more generally, meta-rules.

## THE RULES OF THE MATHEMATICAL GAME

We have seen that students must know what mathematical objects are and their properties, how to represent them and how to use instruments. Is this sufficient to do mathematics? Not at all: using a given instrument to operate on a given representation may be, or not, legitimate (even if done properly). For instance, to solve some numerical problems, some procedures are arithmetical (and are not legitimate in algebra) and other are algebraic (and are not legitimate in arithmetics).

Therefore algebra is not just a question of objects, representations and tools, but also of rules, which are saying what the actions are that we *may* or *may not* do amongst the actions we *can* do. Algebra is not a game in the same sense that chess is a game, but, like chess, algebra does have rules. These rules, moreover, are changing with passing times: the present way of doing differs from, say, the Renaissance Italian way of doing algebra. L. Wittgenstein (the “second Wittgenstein”, the author of the *Philosophical Remarks*, or *On Certainty*, 1986) is an invaluable guide to clarify the extremely complex relationship between objects, signs, practices and rules. (Ernest, 1994, Bagni, 2006).

## SUBPARADIGMS

Some rules (in particular logic) are universal for all mathematics. But other rules are related to a certain domain of mathematics. A square number is always positive, except when studying complex numbers. We call these domains “subparadigms”, which are analogous Kuhn’s paradigms, but less vast, and commensurable between them). This notion of subparadigm allows us to understand the shift from arithmetics to algebra. Semantics (and instrumental value) of the “=” sign change, thus objects (the equalities, the expression with letters) also change. The semiotic systems, although looking quite the same (“2+1 = 3” and “2+$x$ = 3”), are different in fact.

## IDENTIFYING KNOWLEDGE

A last type of knowledge allows us to identify (or recognise) if what we do is mathematical or not, and to identify to what domain of mathematics it belongs. When a student writes something that superficially looks like algebra but actually is wrong or meaningless, the teacher might say: “This is *not* algebra”; and if later the student succeeds in writing a meaningful and correct algebraic text, the teacher might comment: “*This* is algebra”. With these statements, the teacher speaks about the student’s text but also about algebra; he is actually teaching the student *what* is algebra – and what is not[7] (Sackur et al., 2005). We call this Identifying Knowledge; it is also that which allows us to recognise whether a mathematical problem is arithmetical or algebraic, and to choose the appropriate instruments to solve it (without certainty: this kind of knowledge is more abductive that deductive, see Panizza, 2005).

## THE LAYERED DESCRIPTION

As said above, epistemography is not the theory of everything (or, better said, of every kind of knowledge)! Firstly, we only consider here the part of knowledge which is specific to mathematics; this leaves aside nonspecific knowledge, related with the use of (oral and written) natural language or with general reasoning capabilities. “Mathematical activities”, however, remains too vague to allow a precise description. Then, by analogy with the Internet reference model, which is a layered abstract description for the very complex communications and computer network protocol design, we propose a layered description of students’ mathematical activities.

The five descriptive layers of students’ mathematical activities are:

1. the School Layer (what are the students’ rights and duties, why and how to work in the classrooms and at home, what kind of participation is expected by the teacher etc.). This is what french-speaking researchers like Sirota (1993) or Perrenoud (1994) call ““being a student” as a job”[8]. A great number of students’ difficulties may be analysed in terms of the school layer: when they don’t want to learn, or don’t know how to, for instance.

2. The Maths Classroom Layer (how to do maths in the classrooms and at home, what kind of participation is expected by the maths teacher and what is the math teacher supposed to do, etc.). This part of the students’ activities is ruled

[7] which would be almost impossible to do with an explicit discourse within this context: definition or characterization of mathematics are epistemological statements, not mathematical statements

[8] unfortunately, according to Dessus (2004) this concept is almost non existent in English-speaking sociology of education studies.

by what Brousseau (1997) calls the *didactical contract* (see also Sarrazy, 1995, for an extensive survey of this notion). Many students' difficulties can be analysed in terms of didactical contract, as it was brilliantly done by Brousseau (ibid) and followers.

3. The Modelling Layer, which is the description of, for instance, how students change a word problem into a matematical problem, or even how they change a mathematical problem (i. e. expressed in mathematical terms) into an other problem which they can solve with their mathematical tools. A whole field of mathematics education is devoted to the modelling part of the students' mathematical activities (see for instance Lesh and Doerr, 2002).

4. The Discursive Layer, which is the description of students' reasoning on mathematical objects. This reasoning may be expressed by a discourse (like "if $x$ is greater than -3 then $x$+3 is positive and therefore...”), hence the name of this description layer[9]. In France, Duval (2006) is a main contributor in this domain, which is closely related to researches on argumentation (see for instance Yackel and Cobb, 1997) and on proofs (see for instance Gila Hanna, 2000).

5. The deepest, Symbolic Manipulation Layer, describes how students operate on symbolic forms to yield other symbolic forms which represent the solutions of the problem. In the case of algebraic thinking, not too many authors (see for instance Bell, 1996 or Brown & Drouhard, 2003) stress on that – mainly because on the contrary it is often overemphasized by textbooks and teachers.

It is important to notice that what is layered is the description, not the student's activity. It is very similar to what happens in linguistics: the language's description is split in phonetics, syntax, semantics, pragmatics etc. but the subject's act of speech, on the contrary, is of a whole.

## CONCLUSION

A way to cope with the problem of identifying students' mathematical knowledge has long been to focus on students' solving abilities and this can explain the prominent role which has been given to assessment throughout the world. However, many mathematics educators remain reluctant to reduce assessment criteria to solving abilities. Our point is that solving abilities are not so relevant clues on what students know and what they are supposed to know. On the one hand, the student's failure in achieving a task does not give much information on what his or her deficiencies or

9 It is not called "reasoning layer" since that could lead to the erroneous idea that there is no reasoning outside this level.

misconceptions are. On the other hand, the student's success may just show his or her technical abilities, but we cannot be sure that s/he understood conceptually.

Then, how can we determine what students know and are supposed to know? We claim that epistemography can provide accurate answers to this question.